\documentclass[11pt]{amsart}

\usepackage{amsmath}
\usepackage{amssymb}
\usepackage{cite}

\allowdisplaybreaks

\newtheorem{theorem}{Theorem}[section]

\newtheorem{proposition}[theorem]{Proposition}
\newtheorem{remark}[theorem]{Remark}
\theoremstyle{definition}
\newtheorem{definition}[theorem]{Definition}

\newtheorem{notation}[theorem]{Notation}
\newtheorem{assump}[theorem]{Assumption}

\numberwithin{equation}{section}

\AtBeginDocument{{\noindent\small
This is a preprint of a paper whose final and definite form is with 
\emph{Pure and Applied Functional Analysis}, 
ISSN 2189-3756 (print), ISSN 2189-3764 (online). 
Submitted 10-Dec-2017; revised 18-April-2018; accepted for publication 18-April-2018.}
\vspace{9mm}}

\title[A sufficient optimality condition]{A sufficient optimality condition\\ 
for non-linear delayed\\ 
optimal control problems}

\author[A. P. Lemos-Pai\~{a}o]{Ana P. Lemos-Pai\~{a}o}
\address[A. P. Lemos-Pai\~{a}o]{Center for Research and Development 
in Mathematics and Applications (CIDMA),
Department of Mathematics,
University of Aveiro, 
3810-193 Aveiro, Portugal}
\email{{\tt anapaiao@ua.pt}}
 	
\author[C. J. Silva]{Cristiana J. Silva}
\address[C. J. Silva]{Center for Research and Development 
in Mathematics and Applications (CIDMA),
Department of Mathematics,
University of Aveiro, 
3810-193 Aveiro, Portugal}
\email{{\tt cjoaosilva@ua.pt}}

\author[D. F. M. Torres]{Delfim F. M. Torres}
\address[D. F. M. Torres]{Center for Research and Development 
in Mathematics and Applications (CIDMA),
Department of Mathematics,
University of Aveiro, 
3810-193 Aveiro, Portugal}
\email[corresponding author]{{\tt delfim@ua.pt}}

\keywords{Non-linear delayed problems, constant delays in state and control variables, 
sufficient optimality condition, Hamilton--Jacobi equation, optimal control}

\subjclass[2010]{49K15, 49L99}

\begin{document}

\begin{abstract}
We prove a sufficient optimality condition for non-linear optimal control problems 
with delays in both state and control variables. Our result requires 
the verification of a Hamilton--Jacobi partial differential equation 
and is obtained through a transformation that allow us to rewrite 
a delayed optimal control problem as an equivalent non-delayed one.
\end{abstract}

\maketitle

% ----------------------------------------------

\section{Introduction}
\label{section_Introd}

The study of delayed systems, which can be optimized and controlled by a certain control function, 
has a long history and has been developed by many researchers: see, e.g., 
\cite{Bashier,Cacace,Elaiw,Klamka,Mordukhovich1,Mordukhovich2,Mordukhovich3,Mordukhovich4,Stumpf,Xia}. 
Such systems can be called retarded, time-lag, or hereditary processes, and find many applications, 
in diverse fields as biology, chemistry, mechanics, economy and engineering: see, e.g., 
\cite{Bashier,Elaiw,Gollmann2,Ivanov,Klamka,Xu1,Xu2}.

Recent results include Noether type theorems for problems of the calculus 
of variations with time delays \cite{MyID:256,MyID:369,MR3633866},
necessary optimality conditions for quantum \cite{MyID:253}
and Herglotz variational problems with time delays \cite{MyID:319,MyID:342},
as well as delayed optimal control problems with integer 
\cite{MR2745102,MyID:304,MyID:231} and non-integer 
(fractional order) dynamics \cite{MyID:264,MyID:298}.
Applications of such theoretical results are found in biology 
and other natural sciences, e.g., in Tuberculosis \cite{MyID:353}
and HIV \cite{MyID:355,MyID:385}. In the present paper, we establish 
a sufficient optimality condition for an optimal control problem, 
which consists to minimize a cost functional $C[u]$ given by
$$
C[u]=g^0(x(b))+\int_{a}^{b}f^0(t,x(t),x(t-r),u(t),u(t-s))dt
$$
subject to a delayed differential system
$$
\dot{x}(t)=f(t,x(t),x(t-r),u(t),u(t-s))
$$
with given initial functions
\begin{align*}
x(t)&=\varphi(t),\quad t\in[a-r-s,a],\\
u(t)&=\psi(t),\quad t\in[a-s,a[,
\end{align*}
where $r,s>0$, $x(t)\in\mathbb{R}^n$ for each $t\in[a-r-s,b]$, 
$u(t)\in\Omega\subseteq\mathbb{R}^m$ for each $t\in[a-s,b]$ 
and $x(b)\in G\subseteq\mathbb{R}^n$. In order to prove our
sufficient optimality condition, we use a technique 
proposed by Guinn in \cite{Guinn} and used by G\"{o}llmann et al. 
in \cite{Gollmann,Gollmann2}. The technique consists 
to transform a delayed optimal control problem into 
an equivalent non-delayed optimal control problem. After doing 
such transformation, one can apply well-known results for non-delayed 
optimal control problems and then return to the initial delayed problem. 
Analogously to G\"{o}llmann et al. \cite{Gollmann}, we ensure 
the commensurability assumption between the, possibly different, delays 
of state and control variables. We restrict ourselves to delayed 
problems with deterministic controls. For the stochastic case, 
we refer the reader to \cite{Federico2,Fuhrman,Goldys,Ivanov,Larssen}.

Delayed optimal control problems with differential systems, 
which are linear both in state and control, are investigated 
in \cite{Cacace,Chyung,Delfour,Eller,Khellat,Koepcke,Koivo,Lee1,Oguztoreli1,Palanisamy}. 
In some of these papers, necessary and sufficient optimality conditions are derived. 
Our result is different, because we consider a non-linear differential 
system with both delays in state and control variables. Although Banks 
has analyzed non-linear delayed problems, 
he does not consider lags in the control \cite{Banks}. 
Here we consider a delay also in the control variables.

In \cite{Hughes}, Hughes consider variational problems with only 
one constant lag and derives various optimality conditions for them. 
These variational problems can easily be transformed to control problems 
with only one constant delay (see, e.g., \cite[p.~53--54]{Paiao}). 
Hughes also derives an optimality condition for a control problem 
with a constant delay that is the same for the state 
and control variables \cite{Hughes}. Chan and Yung \cite{Chan} 
and Sabbagh \cite{Sabbagh} consider problems that are similar 
to the problems studied by Hughes \cite{Hughes}. In contrast, here 
the state delay is not necessarily equal to the control delay.

In \cite{Jacobs}, Jacobs and Kao investigate delayed problems 
that consist to minimize a cost functional without delays subject 
to a differential system defined by a non-linear function 
with a delay in the state and another one in the control,  
not necessarily equal. Jacobs and Kao begin by transforming their problem 
into a Lagrange-multiplier system subject to a controllability condition
and prove some necessary optimality conditions. Then, they prove existence, 
uniqueness and sufficient conditions in particular situations, 
namely when the differential system is linear in the state 
and control variables \cite{Jacobs}. Here we prove a sufficient 
condition for more general non-linear problems.

As it is well-known, and as Hwang and Bien wrote in \cite{Hwang}, 
many researchers have directed their efforts to seek sufficient 
optimality conditions for control problems with delays: see, e.g., 
\cite{Chyung,Eller,Hughes,Jacobs,Lee3,Schmitendorf}. Therefore,
it is not a surprise that there are authors that already 
proved some sufficient optimality conditions for delayed optimal 
control problems similar but, nevertheless, different from ours. 
In \cite{Schmitendorf}, Schmitendorf consider controls taking values 
in $\mathbb{R}^m$ while here the controls take values in a given set 
$\Omega\subseteq\mathbb{R}^m$, $m\in\mathbb{N}$. Lee and Yung study 
a problem similar to the one considered by Schmitendorf \cite{Schmitendorf}, 
but where the control belongs to a subset of $\mathbb{R}^m$, 
as we do here \cite{Lee3}. However, their sufficient conditions 
are different than our. In particular, \cite{Lee3} assumes 
the existence of a symmetric matrix under some conditions 
that are not easily computable. In \cite{Hwang}, Hwang and Bien 
prove a sufficient condition for problems involving a differential 
affine time-delay system with the same time delay for the state 
and for the control. In 1996, Lee and Yung derived various first 
and second-order sufficient conditions for non-linear optimal control problems 
with only a constant delay in the state \cite{LeeC}. 
Their class of problems is obviously different from our.
In particular, we consider delays for both state 
and control variables. In 2006, Basin and Rodriguez-Gonzalez 
proved a necessary and a sufficient optimality condition 
for a problem that consists to minimize 
a quadratic cost functional subject to a a linear system 
with multiple time delays in the control variable \cite{Basin}. 
In their work, they begin by deriving a necessary condition through 
Pontryagin's maximum principle. Afterwards, sufficiency is proved 
by verifying if the candidate found, through the maximum principle, 
satisfies the Hamilton--Jacobi--Bellman equation. Although Basin 
and Rodriguez-Gonzalez consider multiple time delays, the dependence 
of the state and control in the differential system is linear while
here is non-linear. Later, in 2010 and 2011, Federico et al. 
devoted their attention to optimal control problems that only contain 
delays in the state variables and the dependence on the control is linear 
\cite{Federico1,Federico3}. Also in 2010, Carlier and Tahraqui 
investigated optimal control problems with 
a unique delay in the state \cite{Carlier}. 
The most general results on the area of optimal control 
with delay-differential inclusions in infinite dimensions 
seem those of Mordukhovich et al. 
\cite{Mordukhovich1,Mordukhovich2,Mordukhovich3,Mordukhovich4}. 

The paper is organized as follows. In Section~\ref{sect_preliminaries}, 
we recall a useful sufficient optimality condition for a non-linear 
optimal control problem without delays \cite[p.~347--351]{Lee2}. 
Our result, a sufficient optimality condition
for a non-linear optimal control problem with time lags 
both in state and control variables is then formulated 
in Section~\ref{sect_delayOC} (Theorem~\ref{theo_suf_NLD}).
Its proof is given in Section~\ref{sec:proof}.
We end with Section~\ref{sec:5}, where an example 
that illustrates the obtained theoretical result is given.

% ----------------------------------------------

\section{A non-delayed sufficient optimality condition}
\label{sect_preliminaries}

We recall a well-known sufficient optimality condition
for non-linear optimal control problems without delays. 
Consider the following optimal control problem,
which we denote by $(NL)$:
\begin{equation*}
\min\ \ C[u]=g^0(x(b))+\int_{a}^{b}f^0(t,x(t),u(t))dt
\end{equation*}
subject to the non-linear control system in $\mathbb{R}^n$
\begin{equation}
\label{eq_dif_syst}
\dot{x}(t)=f(t,x(t),u(t))
\end{equation}
with initial boundary condition
\begin{equation}
\label{eq_init_cond}
x(a) = x_a,
\end{equation}
where $x(t)\in\mathbb{R}^n$ and $u(t)\in\Omega\subseteq\mathbb{R}^m$ 
for each $t\in[a,b]$. The functions $f^0$, $f$ and $g^0$ are of class 
$\mathcal{C}^1$ with respect to all its arguments and 
\begin{equation}
\label{eq_end_cond}
x(b)\in G\subseteq\mathbb{R}^n. 
\end{equation}
A pair of functions $(x,u)\in W^{1,\infty}([a,b],\mathbb{R}^n)
\times L^{\infty}([a,b],\mathbb{R}^m)$ that satisfies conditions 
\eqref{eq_dif_syst}--\eqref{eq_end_cond} is said to be
admissible for $(NL)$. 

\begin{notation}
Along all the text, we use the notation $\partial_i f$ to denote the partial 
derivative of a given function $f$ with respect to its $i$th argument. 
For example, $\displaystyle\partial_2f^0(t,x,u)
=\frac{\partial f^0}{\partial x}(t,x,u)$.
\end{notation}

The following theorem provides a sufficient optimality condition for $(NL)$.

\begin{theorem}[See Chapter 5.2, Theorem~7 of \cite{Lee2}] 
\label{theo_suf_nonlinear}
Consider problem $(NL)$. Assume there exists a 
$\mathcal{C}^1(\mathbb{R}^{1+2n})$ feedback control
$u^*(t,x(t),\eta(t,x(t)))$ such that
\begin{equation*}
\begin{split}
\max_{u\in\Omega} H\left(t,x(t),u,\eta(t,x(t))\right)
&=H\left(t,x(t),u^*(t,x(t),\eta(t,x(t))),\eta(t,x(t))\right)\\
&=: H^0\left(t,x(t),\eta(t,x(t))\right)
\end{split}
\end{equation*}
for all $t\in[a,b]$, where
\begin{equation*}
H(t,x,u,\eta)=-f^0(t,x,u) + \eta f(t,x,u).
\end{equation*}
Furthermore, suppose that the $\mathcal{C}^2(\mathbb{R}^{1+n})$
function $S(t,x(t))$, $t\in[a,b]$, is a solution of the 
Hamilton--Jacobi equation
\begin{equation*}
\partial_1 S(t,x(t))+H^0(t,x(t),\partial_2S(t,x(t)))=0
\end{equation*}
with $S(b,x(b))=-g^0(x(b))$, $x(b)\in G$, 
and that the control law
$$
u^*(t,x(t),\partial_2S(t,x(t)))
$$
determines a response $\tilde{x}(t)$ steering $(a,x_a)$ to $(b,G)$. Then, 
$$
\tilde{u}(t)=u^*(t,\tilde{x}(t),\partial_2S(t,\tilde{x}(t)))
$$ 
is an optimal control for $(NL)$ with minimal cost 
$C[\tilde{u}]=-S(a,x_a)$.
\end{theorem}

% -------------------------------------------------------------------

\section{Main Result}
\label{sect_delayOC}

In this paper we are interested in non-linear optimal control problems 
with discrete time delays $r \geq 0$ in the state variable 
$x(t) \in \mathbb{R}^n$ and $s \geq 0$ in the control variable 
$u(t) \in \mathbb{R}^m$, $(r, s) \neq (0, 0)$. Let us define our
problem.

\begin{definition}
\label{def:NLD}
The non-linear delayed optimal control problem $(NLD)$ consists in
\begin{equation*}
\min\ \ C_D[u]=g^0(x(b))+\int_{a}^{b}f^0(t,x(t),x(t-r),u(t),u(t-s))dt
\end{equation*}
subject to the delayed differential system
\begin{equation}
\label{eq:Delay:linear:CS}
\dot{x}(t)=f(t,x(t),x(t-r),u(t),u(t-s))
\end{equation}
with given initial conditions
\begin{align}
x(t)&=\varphi(t),\ t\in[a-r-s,a],\\
\label{eq:InitCond:delayControl}
u(t)&=\psi(t),\ t\in[a-s,a[,
\end{align}
where $x(t)\in\mathbb{R}^n$ for each $t\in[a-r-s,b]$, 
$u(t)\in\Omega\subseteq\mathbb{R}^m$ for each 
$t\in[a-s,b]$ and $x(b)\in G\subseteq\mathbb{R}^n$. 
Functions $f^0$, $f$ and $g^0$ are of class 
$\mathcal{C}^1$ with respect to all their arguments. 
Admissible pairs $(x,u)$ of problem ($NLD$) 
satisfy $(x,u)\in W^{1,\infty}([a-r-s,b],
\mathbb{R}^n)\times L^{\infty}([a-s,b],\mathbb{R}^m)$ 
and all the conditions 
\eqref{eq:Delay:linear:CS}--\eqref{eq:InitCond:delayControl}.
\end{definition}

In what follows, we assume that the time delays 
$r$ and $s$ respect the following assumption.

\begin{assump}[Commensurability assumption]
\label{assump}
We consider $r,s\geq 0$ not simultaneously equal to zero and commensurable, 
that is, 
$$
(r,s)\neq(0,0)
$$ 
and
\begin{equation*}
\frac{r}{s}\in\mathbb{Q}\ \text{ for }\ s>0
\ \text{ or }\ \frac{s}{r}\in\mathbb{Q}\ \text{ for }\ r>0.
\end{equation*}
\end{assump}
Assumption~\ref{assump} holds for any couple of rational numbers $(r, s)$ 
where at least one of them is nonzero \cite{Gollmann}.

\begin{notation}
In what follows,  
$x_a=x(a)=\varphi(a)$;
$x_r(t)=(x(t),x(t-r))$;
$t_s=t-s$; 
and $t^s=t+s$. Moreover, we define the operators
$[\cdot,\cdot]_r$ and $\langle\cdot,\cdot\rangle_r$ by
$[x,\zeta]_r(t) := (t,x_r(t),\zeta(t,x_r(t)))$
and
$\langle x,\zeta\rangle_r(t) := (t,x_r(t),\zeta(t,x(t)))$.
\end{notation}

Our result generalizes Theorem~\ref{theo_suf_nonlinear}
for the non-linear delayed optimal control problem $(NLD)$
of Definition~\ref{def:NLD}.

\begin{theorem}
\label{theo_suf_NLD}
Consider problem $(NLD)$ and let the interval $[a,b]$ 
be divided into $N \in\mathbb{N}$ subintervals of amplitude 
$h = \frac{b-a}{N}>0$. Assume there exists a $\mathcal{C}^1(\mathbb{R}^{1+3n})$ 
feedback control $u^*(t,x_r(t),\eta(t,x_r(t))) = u^*[x,\eta]_r(t)$ such that
\begin{equation}
\label{eq_max_cond_main_theo}
\begin{split}
&\max_{u\in\Omega} \bigl\{H(t,x_r(t),u,u^*[x,\eta]_r(t_s),\eta(t,x_r(t)))\\
&\quad + H(t^s,x_r(t^s),u^*[x,\eta]_r(t^s),
u,\eta(t^s,x_r(t^s)))\chi_{[a,b-s]}(t)\bigr\}\\
&= H(t,x_r(t),u^*[x,\eta]_r(t),u^*[x,\eta]_r(t_s),\eta(t,x_r(t)))\\
&\quad +H(t^s,x_r(t^s),u^*[x,\eta]_r(t^s),
u^*[x,\eta]_r(t),\eta(t^s,x_r(t^s)))\chi_{[a,b-s]}(t)\\
&=: H^0[x,\eta]_r(t)+H^0[x,\eta]_r(t^s)\chi_{[a,b-s]}(t)
\end{split}
\end{equation}
for all $t\in[a,b]$, where
\begin{equation*}
H(t,x,y,u,v,\eta)=-f^0(t,x,y,u,v)+ \eta f(t,x,y,u,v).
\end{equation*}
Furthermore, let $I_i=[a+hi,a+h(i+1)]$, 
$i=0, \ldots, N-1$, and suppose that function 
$S(t,x(t))\in\mathcal{C}^2(\mathbb{R}^{1+n})$, $t\in[a,b]$, 
is a solution of equation
\begin{multline}
\label{eq_HJ_main_theo}
\partial_1S(t,x(t))+\sum_{i=0}^{N-1}\big\{-f^0(t,x_r(t),u^*\langle x,\partial_2 S\rangle_r(t),
u^*\langle x,\partial_2 S\rangle_r(t_s))\\
+\partial_2S(t,x(t))f(t,x_r(t),u^*\langle x,\partial_2 S\rangle_r(t),
u^*\langle x,\partial_2 S\rangle_r(t_s))\big\}\chi_{I_i}(t)=0
\end{multline}
with $S(b,x(b))=-g^0(x(b))$, $x(b)\in G$.
Finally, consider that the control law
$$
u^*\left(t,x_r(t),\partial_2S(t,x(t))\right) 
= u^*\langle x,\partial_2 S\rangle_r(t),
\quad t\in[a,b],
$$
determines a response $\tilde{x}(t)$ 
steering $(a,x_a)$ to $(b,G)$. Then, 
$$
\tilde{u}(t)
=u^*\left(t,\tilde{x}(t),\tilde{x}(t-r),\partial_2S(t,\tilde{x}(t)\right)
$$
is an optimal control for $(NLD)$ that leads to the minimal cost 
$$
C_D[\tilde{u}]=-S(a,x_a).
$$
\end{theorem}

We prove Theorem~\ref{theo_suf_NLD} in Section~\ref{sec:proof}.

% ---------------------------------------

\section{Proof of the delayed sufficient optimality condition}
\label{sec:proof}

We prove Theorem~\ref{theo_suf_NLD} as a corollary of
Theorem~\ref{theo_suf_nonlinear} by
transforming the non-linear delayed optimal control problem 
$(NLD)$ into an equivalent non-linear optimal control problem 
without delays of type $(NL)$. For that, 
we use the approach of \cite{Gollmann,Guinn}.
Without loss of generality, we assume the first case of Assumption~\ref{assump}, 
that is, $\displaystyle\frac{r}{s}\in\mathbb{Q}$ for $r\geq0$ and $s>0$.
Consequently, there exist $k,l\in\mathbb{N}$ such that
\begin{equation*}
\frac{r}{s}=\frac{k}{l}
\Leftrightarrow
rl=sk
\Leftrightarrow
\frac{r}{k}=\frac{s}{l}.
\end{equation*}
Thus, we divide the interval $[a,b]$ into $N\in\mathbb{N}$ 
subintervals of amplitude $h:=\displaystyle\frac{r}{k}=\frac{s}{l}$. 
We can note that 
$$
r=hk \text{ and } s=hl.
$$
Furthermore, let us assume that 
\begin{equation*}
a+hN=b\ \text{ and }\ N>2k+1,
\end{equation*}
with $N\in\mathbb{N}$.

\begin{remark}
If $b-a$ is not a multiple of $h$, that is, $b-a\neq hN$, 
then we can study problem ($NLD$) for $t\in[a,\tilde{b}]$, 
where $\tilde{b}$ is the smallest multiple of $h$ 
greater than $b$. Thus, we also study problem ($NLD$) 
for $t\in[a,b]$, because $b<\tilde{b}$.
\end{remark}

For $i=0,\ldots,N-1$ and $t\in[a,a+h]$, we define the new variables
\begin{align*}
\xi_i(t)=x(t+hi)\ \text{ and }\
\theta_i(t)=u(t+hi).
\end{align*}
The non-linear delayed problem $(NLD)$ is transformed into the 
following equivalent non-linear problem $(\overline{NL})$ 
without delays: 
\begin{equation}
\label{cost_function_NL-}
\min \overline{C}[\theta]=g^0(\xi_{N-1}(a+h))
+\int_{a}^{a+h}\sum_{i=0}^{N-1}f^0(t+hi,\xi_i(t),\xi_{i-k}(t),\theta_i(t),\theta_{i-l}(t))dt
\end{equation}
subject to the non-delayed differential system
\begin{equation}
\label{diff_syst_NL-}
\dot{\xi_i}(t)=f(t+hi,\xi_i(t),\xi_{i-k}(t),\theta_i(t),\theta_{i-l}(t)),
\ i=0,\ldots,N-1, \  t\in[a,a+h],
\end{equation}
and the initial conditions
\begin{equation}
\label{eq_continuity_cond_NL-}
\begin{split}
\xi_i(t)&=\varphi(t+hi),\ i=-k-l,\ldots,-1, \ t\in[a,a+h],\\
\theta_i(t)&=\psi(t+hi),\ i=-l,\ldots,-1, \ t\in[a,a+h[,\\
\xi_i(a+h)&=\xi_{i+1}(a),\ i=0,\ldots,N-2.
\end{split}
\end{equation}
We observe that the cost functional \eqref{cost_function_NL-} 
depends only on $t\in[a,a+h]$, 
$\xi(t)=[\xi_0(t)\ \ldots\ \xi_{N-1}(t)]^T$ and 
$\theta(t)=[\theta_0(t)\ \ldots\ \theta_{N-1}(t)]^T$, because 
$$
\xi^-(t)=[\xi_{-k-l}(t)\ \xi_{1-k-l}(t)\ \ldots\ \xi_{-1}(t)]^T
$$ 
and 
$$
\theta^-(t)=[\theta_{-l}(t)\ \ldots\ \theta_{-1}(t)]^T
$$ 
are already known. Thus, the integrand function of 
\eqref{cost_function_NL-} can be written as
\begin{equation*}
\sum_{i=0}^{N-1}f^0\left(t+hi,\xi_i(t),\xi_{i-k}(t),
\theta_i(t),\theta_{i-l}(t)\right)
=F^0\left(t,\xi(t),\theta(t)\right).
\end{equation*}
We can also write
\begin{equation*}
g^0(\xi_{N-1}(a+h))=G^0(\xi(a+h)).
\end{equation*}
Note that we are writing $G^0$ as a function of $\xi(a+h)\in\mathbb{R}^{nN}$ 
in order to obtain problem $(\overline{NL})$ written in the form 
used by Theorem~\ref{theo_suf_nonlinear}. However, function $G^0$ 
depends only on $\xi_{N-1}(a+h)\in\mathbb{R}^n$.
Consequently, we have
\begin{multline*}
g^0(\xi_{N-1}(a+h))+\int_{a}^{a+h}\sum_{i=0}^{N-1}
f^0(t+hi,\xi_i(t),\xi_{i-k}(t),\theta_i(t),\theta_{i-l}(t))dt\\
=G^0(\xi(a+h))+\int_{a}^{a+h}F^0(t,\xi(t),\theta(t))dt.
\end{multline*}
Using similar arguments, the differential system 
\eqref{diff_syst_NL-} can be written as
\begin{align*}
\dot{\xi}(t) &= 
\left[\begin{matrix}
\dot{\xi}_0(t)\\
\dot{\xi}_1(t)\\
\vdots\\
\dot{\xi}_{N-1}(t)\\
\end{matrix}\right]\\
&=
\left[\begin{matrix}
f(t,\xi_0(t),\xi_{-k}(t),\theta_0(t),\theta_{-l}(t))\\
f(t+h,\xi_1(t),\xi_{1-k}(t),\theta_1(t),\theta_{1-l}(t))\\
\vdots\\
f\left(t+h(N-1),\xi_{N-1}(t),\xi_{N-1-k}(t),\theta_{N-1}(t),\theta_{N-1-l}(t)\right)\\
\end{matrix}\right]\\
&= F(t,\xi(t),\theta(t)), \quad t\in[a,a+h].
\end{align*}
In order to apply Theorem~\ref{theo_suf_nonlinear}, 
we consider the initial boundary condition, 
with respect to variable $\xi$, given by
\begin{equation*}
\xi_a=\xi(a)=
\left[\begin{matrix}
\xi_0(a)\\
\xi_1(a)\\
\vdots\\
\xi_{N-1}(a)
\end{matrix}\right]
=
\left[\begin{matrix}
x_a\\
\xi_0(a+h)\\
\vdots\\
\xi_{N-2}(a+h)
\end{matrix}\right].		
\end{equation*}
\begin{remark}
Only the first component of $\xi_a$ is known a priori. The others 
are determined using the continuity conditions 
$\xi_i(a+h)=\xi_{i+1}(a)$ of \eqref{eq_continuity_cond_NL-}, 
$i=0,\ldots,N-2$, and the fixed value $x_a$. 
\end{remark}
Concluding, problem $(\overline{NL})$ is written 
in the standard form, as follows:
\begin{align*}
\min\ \ &\overline{C}[\theta]=G^0(\xi(a+h))+\int_{a}^{a+h}F^0(t,\xi(t),\theta(t))dt\\
\text{s.t. }& \dot{\xi}(t)=F(t,\xi(t),\theta(t)),\ t\in[a,a+h],\\
&\xi(a) = \xi_a =
\left[\begin{matrix}
x_a\\
\xi_0(a+h)\\
\vdots\\
\xi_{N-2}(a+h)
\end{matrix}\right],
\end{align*}
knowing $\xi^-(t)$ and $\theta^-(t)$ for all $t\in[a,a+h]$ 
and ensuring the continuity conditions
$\xi_i(a+h)=\xi_{i+1}(a)$ of \eqref{eq_continuity_cond_NL-}, 
$i=0,\ldots,N-2$. Furthermore, we know that
\begin{itemize}
\item $\xi(t)\in\mathbb{R}^{nN}$ and $\theta(t)\in\overline{\Omega}
\subseteq\mathbb{R}^{mN}$ for each $t\in[a,a+h]$;
\item $\xi(a+h)\in \overline{G}=\mathbb{R}^{n(N-1)}\times G$;
\item functions $F^0$, $F$ and $G^0$ are of class $\mathcal{C}^1$ 
with respect to all their arguments, because $f^0$, 
$f$ and $g^0$ are of class $\mathcal{C}^1$ in all their arguments.
\end{itemize}
Therefore, we are in condition to apply Theorem~\ref{theo_suf_nonlinear}. 
Firstly, we are going to prove the first part of Theorem~\ref{theo_suf_NLD}, 
that is, we show that \eqref{eq_max_cond_main_theo} holds. 
Assume there exists a feedback control 
$\theta^*(t,\xi(t),\Lambda(t,\xi(t)))\in\mathcal{C}^1(\mathbb{R}^{1+2nN})$ 
such that
\begin{equation}
\label{eq_maximality_NL-}
\begin{split}
\max_{\theta\in\overline{\Omega}}\overline{H}(t,\xi(t),\theta,\Lambda(t,\xi(t)))
&=\overline{H}(t,\xi(t),\theta^*(t,\xi(t),\Lambda(t,\xi(t))),\Lambda(t,\xi(t)))\\
&=:\overline{H}^0(t,\xi(t),\Lambda(t,\xi(t)))
\end{split}
\end{equation}
for all $t\in[a,a+h]$, where $\overline{H}(t,\xi,\theta,\Lambda)
=-F^0(t,\xi,\theta)+ \Lambda F(t,\xi,\theta)$. 
In order to write the previous 
condition with respect to the original variables, 
we do the following remark.

\begin{remark}
For each $t\in[a,b]$, there exists $j\in\{0,\ldots,N-1\}$ such that 
$$
a+hj\leq t \leq a+h(j+1) \Leftrightarrow a\leq t-hj \leq a+h.
$$ 
Thus, let us define $t'\in[a,a+h]$ as $t'=t-hj$ 
and $\eta(t, x(t), x(t-r))$ as
\begin{equation*}
\eta(t, x(t+hj), x(t+hj-r))=\Lambda^j(t-hj, \xi_j(t), \xi_{j-k}(t)).
\end{equation*}
Then, 
\begin{align*}
\Lambda^j(t,\xi_j(t),\xi_{j-k}(t))
&=\Lambda^j(t+hj-hj,\xi_j(t),\xi_{j-k}(t))\\
&=\eta(t+hj,x(t+hj),x(t+hj-r)),
\end{align*}
that is,
\begin{equation*}
\begin{split}
\Lambda^j(t',\xi_j(t'),\xi_{j-k}(t'))
&=\eta(t'+hj,x(t'+hj),x(t'+hj-r))\\
&=\eta(t,x(t),x(t-r))
\end{split}
\end{equation*}
and 
\begin{align*}
\Lambda^{j+l}&(t,\xi_{j+l}(t),\xi_{j+l-k}(t))\\
&= \Lambda^{j+l}(t+h(j+l)-h(j+l),x(t+hj+hl),x(t+hj+hl-hk))\\
&= \Lambda^{j+l}(t+hj+s-h(j+l),x(t+hj+s),x(t+hj+s-r))\\
&= \eta(t+hj+s,x(t+hj+s),x(t+hj+s-r)),
\end{align*}
which implies that
\begin{align*}
\Lambda^{j+l}&\left(t',\xi_{j+l}(t'),\xi_{j+l-k}(t')\right)\\
&=\eta\left(t'+hj+s,x(t'+hj+s),x(t'+hj+s-r)\right)\\
&=\eta\left(t+s,x(t+s),x(t+s-r)\right).
\end{align*}	
\end{remark} 

As equation~\eqref{eq_maximality_NL-} is verified for all admissible 
$\theta\in\overline{\Omega}$, we can choose an admissible control 
$\overline{\theta}\in\overline{\Omega}$ such that
\begin{align}
\label{admissible_control}
\overline{\theta}_i=
\begin{cases}
\theta^*_i(t',\xi_i(t'),\xi_{i-k}(t'),\Lambda^i(t',\xi_i(t'),\xi_{i-k}(t'))), & i\neq j,\\
\theta_i, & i=j,
\end{cases}
\end{align}
$i=0,\ldots,N-1$, where 
$\theta=[
\begin{matrix}
\theta_0 & \ldots & \theta_{N-1}
\end{matrix}]^T$ 
is an admissible control for problem $(\overline{NL})$. 
From condition \eqref{eq_maximality_NL-}, we can write that
\begin{equation*}
\overline{H}(t',\xi(t'),\overline{\theta},\Lambda(t',\xi(t')))
\leq\overline{H}\left(t',\xi(t'),\theta^*(t',\xi(t'),
\Lambda(t',\xi(t'))),\Lambda(t',\xi(t'))\right).
\end{equation*}	
From now on, we write $\theta^{*'}$ instead of 
$\theta^*(t',\xi(t'),\Lambda(t',\xi(t')))$, in order 
to simplify expressions. With this notation, we have
\begin{multline*}
-F^0\left(t',\xi(t'),\overline{\theta}\right)
+\Lambda\left(t',\xi(t')\right) 
F\left(t',\xi(t'),\overline{\theta}\right)\\
\leq -F^0\left(t',\xi(t'),\theta^{*'}\right)
+\Lambda\left(t',\xi(t')\right) 
F\left(t',\xi(t'),\theta^{*'}\right),
\end{multline*}
which is equivalent to
\begin{align*}
\sum_{i=0}^{N-1}&\Bigl\{-f^0\left(t'+hi,\xi_i(t'),\xi_{i-k}(t'),
\overline{\theta}_i,\overline{\theta}_{i-l}\right)\\
&\quad +\Lambda^i\left(t',\xi_i(t'),\xi_{i-k}(t')\right)
f\left(t'+hi,\xi_i(t'),\xi_{i-k}(t'),\overline{\theta}_i,\overline{\theta}_{i-l}\right)\Bigr\}\\
&\leq \sum_{i=0}^{N-1}\Bigl\{-f^0\left(t'+hi,\xi_i(t'),
\xi_{i-k}(t'),\theta^{*'}_i,\theta^{*'}_{i-l}\right)\\
&\quad +\Lambda^i\left(t',\xi_i(t'),\xi_{i-k}(t')\right)
f\left(t'+hi,\xi_i(t'),\xi_{i-k}(t'),\theta^{*'}_i,\theta^{*'}_{i-l}\right)\Bigr\}.
\end{align*}
Considering $I=\{0,\ldots,N-1\}\backslash\{j,j+l\}$ 
and definition \eqref{admissible_control} for 
the admissible control $\overline{\theta}$, we obtain that
\begin{align*}
\sum_{i\in I}&\Bigl\{-f^0\left(t'+hi,\xi_i(t'),\xi_{i-k}(t'),\theta^{*'}_i,\theta^{*'}_{i-l}\right)\\
&\quad +\Lambda^i\left(t',\xi_i(t'),\xi_{i-k}(t')\right)f\left(t'+hi,\xi_i(t'),
\xi_{i-k}(t'),\theta^{*'}_i,\theta^{*'}_{i-l}\right)\Bigr\}\\
&\quad -f^0\left(t'+hj,\xi_j(t'),\xi_{j-k}(t'),\theta_j,\theta^{*'}_{j-l}\right)\\
&\quad +\Lambda^j\left(t',\xi_j(t'),\xi_{j-k}(t')\right)
f\left(t'+hj,\xi_j(t'),\xi_{j-k}(t'),\theta_j,\theta^{*'}_{j-l}\right)\\
&\quad + \Bigl[-f^0\left(t'+hj+s,\xi_{j+l}(t'),\xi_{j+l-k}(t'),\theta^{*'}_{j+l},\theta_j\right)\\
&\qquad +\Lambda^{j+l}\left(t',\xi_{j+l}(t'),\xi_{j+l-k}(t')\right)\\
&\qquad \times f\left(t'+hj+s,\xi_{j+l}(t'),
\xi_{j+l-k}(t'),\theta^{*'}_{j+l},\theta_j\right)\Bigr]\chi_{\{0,\ldots,N-1-l\}}(j)\\
&\leq \sum_{i=0}^{N-1}\Bigl\{-f^0\left(t'+hi,\xi_i(t'),
\xi_{i-k}(t'),\theta^{*'}_i,\theta^{*'}_{i-l}\right)\\
&\quad +\Lambda^i\left(t',\xi_i(t'),\xi_{i-k}(t')\right)
f\left(t'+hi,\xi_i(t'),\xi_{i-k}(t'),\theta^{*'}_i,\theta^{*'}_{i-l}\right)\Bigr\}.
\end{align*}
The terms of the first and second members with indexes 
in set $I$ cancel, and we simply have
\begin{equation} 
\label{eq_maximality_cond_desenvolvida}
\begin{split}
&-f^0\left(t'+hj,\xi_j(t'),\xi_{j-k}(t'),\theta_j,\theta^{*'}_{j-l}\right)\\
&\quad +\Lambda^j\left(t',\xi_j(t'),\xi_{j-k}(t')\right)\\
&\quad \times f\left(t'+hj,\xi_j(t'),\xi_{j-k}(t'),\theta_j,\theta^{*'}_{j-l}\right)\\
&\quad + \Bigl[-f^0\left(t'+hj+s,\xi_{j+l}(t'),\xi_{j+l-k}(t'),\theta^{*'}_{j+l},\theta_j\right)\\
&\qquad+\Lambda^{j+l}\left(t',\xi_{j+l}(t'),\xi_{j+l-k}(t')\right)\\
&\qquad \times f\left(t'+hj+s,\xi_{j+l}(t'),\xi_{j+l-k}(t'),\theta^{*'}_{j+l},\theta_j\right)\Bigr]
\chi_{\{0,\ldots,N-1-l\}}(j)\\
&\leq -f^0\left(t'+hj,\xi_j(t'),\xi_{j-k}(t'),\theta^{*'}_j,\theta^{*'}_{j-l}\right)\\
&\quad +\Lambda^j\left(t',\xi_j(t'),\xi_{j-k}(t')\right)
f(t'+hj,\xi_j(t'),\xi_{j-k}(t'),\theta^{*'}_j,\theta^{*'}_{j-l})\\
&\quad + \Bigl[-f^0\left(t'+hj+s,\xi_{j+l}(t'),\xi_{j+l-k}(t'),\theta^{*'}_{j+l},\theta^{*'}_j\right)\\
&\qquad +\Lambda^{j+l}\left(t',\xi_{j+l}(t'),\xi_{j+l-k}(t')\right)\\
&\qquad \times f\left(t'+hj+s,\xi_{j+l}(t'),\xi_{j+l-k}(t'),\theta^{*'}_{j+l},
\theta^{*'}_j\right)\Bigr]\chi_{\{0,\ldots,N-1-l\}}(j).
\end{split}
\end{equation}
We can observe that
\begin{equation*}
\begin{gathered}
t'+hj=t-hj+hj=t;\\
\xi_j(t')=x(t'+hj)=x(t);\\
\xi_{j-k}(t')=x(t'+hj-hk)=x(t-r);\\
\xi_{j-l}(t')=x(t'+hj-hl)=x(t-s);\\
\xi_{j-l-k}(t')=x(t'+hj-hl-hk)=x(t-s-r);\\
\xi_{j+l}(t')=x(t'+hj+hl)=x(t+s);\\
\xi_{j+l-k}(t')=x(t'+hj+hl-hk)=x(t+s-r);
\end{gathered}
\end{equation*}
\begin{align*}
\theta^{*'}_j
&= \theta^*_j\left(t',\xi_j(t'),\xi_{j-k}(t'),
\Lambda^j\left(t',\xi_j(t'),\xi_{j-k}(t')\right)\right)\\
&= u^*\left(t'+hj,x(t),x(t-r),\eta\left(t,x(t),x(t-r)\right)\right)\\
&= u^*[x,\eta]_r(t);\\
\theta^{*'}_{j-l}
&= \theta^*_{j-l}\left(t',\xi_{j-l}(t'),\xi_{j-l-k}(t'),
\Lambda^{j-l}\left(t',\xi_{j-l}(t'),\xi_{j-l-k}(t')\right)\right)\\
&= u^*\left(t'+hj-hl,x(t-s),x(t-s-r),\eta(t-s,x(t-s),x(t-s-r))\right)\\
&= u^*\left(t_s,x_r(t_s),\eta(t_s,x_r(t_s))\right)= u^*[x,\eta]_r(t_s);
\end{align*}
\begin{align*}
\theta^{*'}_{j+l}
&= \theta^*_{j+l}\left(t',\xi_{j+l}(t'),\xi_{j+l-k}(t'),
\Lambda^{j+l}\left(t',\xi_{j+l}(t'),\xi_{j+l-k}(t')\right)\right)\\
&= u^*\left(t'+hj+hl,x(t+s),x(t+s-r),\eta(t+s,x(t+s),x(t+s-r))\right)\\
&= u^*\left(t^s,x_r(t^s),\eta(t^s,x_r(t^s))\right) = u^*[x,\eta]_r(t^s);
\end{align*}
$\theta_j=u$, where $u\in\Omega$ 
is an arbitrary admissible control of problem $(NLD)$.	  
Using these relations, we rewrite the first member of inequality
\eqref{eq_maximality_cond_desenvolvida} as
\begin{align*}
-f^0&\left(t'+hj,\xi_j(t'),\xi_{j-k}(t'),\theta_j,\theta^{*'}_{j-l}\right)\\
&\quad +\Lambda^j\left(t',\xi_j(t'),\xi_{j-k}(t')\right)
f\left(t'+hj,\xi_j(t'),\xi_{j-k}(t'),\theta_j,\theta^{*'}_{j-l}\right)\\
&\quad +\Bigl[-f^0(t'+hj+s,\xi_{j+l}(t'),\xi_{j+l-k}(t'),\theta^{*'}_{j+l},\theta_j)\\
&\qquad +\Lambda^{j+l}\left(t',\xi_{j+l}(t'),\xi_{j+l-k}(t')\right)\\
&\qquad \times f\left(t'+hj+s,\xi_{j+l}(t'),
\xi_{j+l-k}(t'),\theta^{*'}_{j+l},\theta_j\right)\Bigr]\chi_{\{0,\ldots,N-1-l\}}(j)\\
&= -f^0\left(t,x_r(t),u,u^*[x,\eta]_r(t_s)\right)\\
&\quad +\eta\left(t,x_r(t)\right)f\left(t,x_r(t),u,u^*[x,\eta]_r(t_s)\right)\\
&\quad + \Bigl[-f^0\left(t^s,x_r(t^s),u^*[x,\eta]_r(t^s),u\right)\\
&\qquad +\eta\left(t^s,x_r(t^s)\right)
f\left(t^s,x_r(t^s),u^*[x,\eta]_r(t^s),u\right)\Bigr]\chi_{[a,b-s]}(t)\\
&= H\left(t,x_r(t),u,u^*[x,\eta]_r(t_s),\eta\left(t,x_r(t)\right)\right)\\
&\quad +H\left(t^s,x_r(t^s),u^*[x,\eta]_r(t^s),u,\eta(t^s,x_r(t^s))\right)
\chi_{[a,b-s]}(t).
\end{align*}
On the other hand, the second member of inequality
\eqref{eq_maximality_cond_desenvolvida} takes the form
\begin{align*}
&-f^0\left(t'+hj,\xi_j(t'),\xi_{j-k}(t'),\theta^{*'}_j,\theta^{*'}_{j-l}\right)\\
&\quad +\Lambda^j\left(t',\xi_j(t'),\xi_{j-k}(t')\right)
f\left(t'+hj,\xi_j(t'),\xi_{j-k}(t'),\theta^{*'}_j,\theta^{*'}_{j-l}\right)\\
&\quad +\Bigl[-f^0\left(t'+hj+s,\xi_{j+l}(t'),\xi_{j+l-k}(t'),\theta^{*'}_{j+l},\theta^{*'}_j\right)\\
&\qquad +\Lambda^{j+l}\left(t',\xi_{j+l}(t'),\xi_{j+l-k}(t')\right)\\
&\qquad \times f\left(t'+hj+s,\xi_{j+l}(t'),\xi_{j+l-k}(t'),\theta^{*'}_{j+l},\theta^{*'}_j\right)\Bigr]
\chi_{\{0,\ldots,N-1-l\}}(j)\\
&=-f^0\left(t,x_r(t),u^*[x,\eta]_r(t),u^*[x,\eta]_r(t_s)\right)\\
&\quad +\eta\left(t,x_r(t)\right)
f\left(t,x_r(t),u^*[x,\eta]_r(t),u^*[x,\eta]_r(t_s)\right)\\
&\quad +\Bigl[-f^0\left(t^s,x_r(t^s),u^*[x,\eta]_r(t^s),u^*[x,\eta]_r(t)\right)\\
&\qquad +\eta\left(t^s,x_r(t^s)\right)
f\left(t^s,x_r(t^s),u^*[x,\eta]_r(t^s),u^*[x,\eta]_r(t)\right)\Bigr]\chi_{[a,b-s]}(t)\\
&= H\left(t,x_r(t),u^*[x,\eta]_r(t),u^*[x,\eta]_r(t_s),\eta(t,x_r(t))\right)\\
&\quad +H\left(t^s,x_r(t^s),u^*[x,\eta]_r(t^s),u^*[x,\eta]_r(t),\eta(t^s,x_r(t^s))\right)\chi_{[a,b-s]}(t).
\end{align*}
Therefore, the inequality~\eqref{eq_maximality_cond_desenvolvida} is equivalent to
\begin{align*}
H&\left(t,x_r(t),u,u^*[x,\eta]_r(t_s),
\eta\left(t,x_r(t)\right)\right)\\
&\quad +H\left(t^s,x_r(t^s),u^*[x,\eta]_r(t^s),u,\eta\left(t^s,x_r(t^s)\right)\right)\chi_{[a,b-s]}(t)\\
& \leq H\left(t,x_r(t),u^*[x,\eta]_r(t),u^*[x,\eta]_r(t_s),\eta(t,x_r(t))\right)\\
&\quad +H\left(t^s,x_r(t^s),u^*[x,\eta]_r(t^s),u^*[x,\eta]_r(t),
\eta\left(t^s,x_r(t^s)\right)\right)\chi_{[a,b-s]}(t),	
\end{align*}
where $u\in\Omega$ is an arbitrary admissible control of problem $(NLD)$. 
We just proved condition \eqref{eq_max_cond_main_theo}. Now we proceed
by proving equation~\eqref{eq_HJ_main_theo}. Let us suppose that function 
$\overline{S}(t,\xi(t))\in\mathcal{C}^2(\mathbb{R}^{1+nN})$, $t\in[a,a+h]$, 
is a solution to the Hamilton--Jacobi equation
\begin{equation}
\label{eq:HJeq:Proof}
\partial_1 \overline{S}(t,\xi(t))
+\overline{H}^0(t,\xi(t),\partial_2\overline{S}(t,\xi(t)))=0
\end{equation}
with $\overline{S}(a+h,\xi(a+h))=-G^0(\xi(a+h))$ for $\xi(a+h)\in \overline{G}$. 
Now, in order to simplify the notation, we write 
\begin{itemize}
\item $\theta^*$ instead of $\theta^*\left(t,\xi(t),\partial_2\overline{S}(t,\xi(t))\right)$;
\item $\theta^*_i$ instead of $\theta^*_i\left(t,\xi_i(t),\xi_{i-k}(t),
\partial_{i+2}\overline{S}(t,\xi(t))\right)$, for $i=0,\ldots,N-1$.
\end{itemize}
Therefore, the Hamilton--Jacobi equation \eqref{eq:HJeq:Proof} is equivalent to
\begin{align*}
&\partial_1 \overline{S}(t,\xi(t))
+\overline{H}\left(t,\xi(t),\theta^*,\partial_2\overline{S}(t,\xi(t))\right)=0\\
\Leftrightarrow&\partial_1 \overline{S}(t,\xi(t))-F^0\left(t,\xi(t),\theta^*\right)
+\partial_2\overline{S}\left(t,\xi(t)\right)
F\left(t,\xi(t),\theta^*\right)=0\\
\Leftrightarrow&\partial_1\overline{S}(t,\xi(t))
+\sum_{i=0}^{N-1}\Bigl\{-f^0\left(t+hi,\xi_i(t),\xi_{i-k}(t),\theta^*_i,\theta^*_{i-l}\right)\\
&+\partial_{i+2}\overline{S}(t,\xi(t))
f\left(t+hi,\xi_i(t),\xi_{i-k}(t),\theta^*_i,\theta^*_{i-l}\right)\Bigr\}=0.
\end{align*}
For all $t\in[a,a+h]$, one has
\begin{equation*}
\begin{split}
\overline{S}(t,\xi(t))
&=\overline{S}\left(t,\xi_0(t),\xi_1(t),\ldots,\xi_{N-1}(t)\right)\\
&=\overline{S}\left(t,x(t),x(t+h),\ldots,x(t+hN-h)\right).
\end{split}
\end{equation*}
So, we can simply write $\overline{S}(t,\xi(t))$ for $t\in[a,a+h]$ 
as a function of $t$ and $x(t)$ for all $t\in[a,b]$:
\begin{align*}
\overline{S}(t,\xi(t))|_{t\in[a,a+h]}:=S(t,x(t))|_{t\in[a,b]}.
\end{align*}
We can also observe that
\begin{align*}
\partial_{i+2}\overline{S}(t,\xi(t))=\partial_2S(t,x(t))\chi_{I_i}(t),
\end{align*}
\begin{multline*} 
f^0\left(t+hi,\xi_i(t),\xi_{i-k}(t),\theta^*_i,\theta^*_{i-l}\right)\\
= f^0\left(t,x_r(t),u^*\langle x,\partial_2 S\rangle_r(t),
u^*\langle x,\partial_2 S\rangle_r(t_s)\right)\chi_{I_i}(t),
\end{multline*}
and
\begin{multline*} 
f\left(t+hi,\xi_i(t),\xi_{i-k}(t),\theta^*_i,\theta^*_{i-l}\right)\\
= f\left(t,x_r(t),u^*\langle x,\partial_2 S\rangle_r(t),
u^*\langle x,\partial_2 S\rangle_r(t_s)\right)\chi_{I_i}(t),
\end{multline*}
for $i=0,\ldots,N-1$. Therefore, we obtain
\begin{align*}
&\partial_1S(t,x(t))+\sum_{i=0}^{N-1}\Bigl\{-f^0\left(t,x_r(t),
u^*\langle x,\partial_2S\rangle_r(t),u^*\langle x,\partial_2 S\rangle_r(t_s)\right)\\
&+\partial_2S(t,x(t))f\left(t,x_r(t),u^*\langle x,\partial_2 S\rangle_r(t),
u^*\langle x,\partial_2 S\rangle_r(t_s)\right)\Bigr\}\chi_{I_i}(t)=0.
\end{align*}
Furthermore, we have to ensure that
\begin{align*}
&\overline{S}(a+h,\xi(a+h))=-G^0(\xi(a+h))\\
\Leftrightarrow & \overline{S}\left(a+h,x(a+h), x(a+2h),\ldots, x(b)\right)
=-g^0\left(\xi_{N-1}(a+h)\right)\\
\Leftrightarrow & 
\overline{S}\left(a+h,x(a+h), x(a+2h),\ldots, x(b)\right)
=-g^0(x(b)),
\end{align*}
which implies that
$$
S(b,x(b))=-g^0(x(b)).
$$
As $\xi(a+h)\in\overline{G}=\mathbb{R}^{n(N-1)}\times G$, 
then $\xi_{N-1}(a+h)=x(b)\in G$. Therefore, we obtain
equation~\eqref{eq_HJ_main_theo} and its conditions 
$S(b,x(b))=-g^0(x(b))$, $x(b)\in G$.
To finish the proof, let us assume that the control law
$$
\theta_i^*\left(t,\xi_i(t),\xi_{i-k}(t),
\partial_{i+2}\overline{S}(t,\xi(t))\right)
=u^*\langle x,\partial_2 S\rangle_r(t)\chi_{I_i}(t)
$$
determines a response $\tilde{\xi}(t)$, $t\in[a,a+h]$, 
steering $(a,\xi_i(a))$ to $(a+h,\overline{G})$, $i=0,\ldots,N-1$. 
Such assumption implies that the control law
$u^*\langle x,\partial_2 S\rangle_r(t)$
determines a response $\tilde{x}(t)$ steering $(a,x_a)$ to $(b,G)$, 
for all $t\in[a,b]$. For $i=0,\ldots,N-1$ and $t\in[a,a+h]$,
$$
\tilde{\theta}_i(t)=\theta_i^*\left(t,\tilde{\xi}_i(t),
\tilde{\xi}_{i-k}(t),\partial_{i+2}\overline{S}(t,\tilde{\xi}(t))\right)
$$
is the $i$th component of an optimal control $\tilde{\theta}(t)$ 
that lead us to the minimal cost
$$
\overline{C}[\tilde{\theta}]=-\overline{S}(a,\xi(a))
=-\overline{S}\left(a,\xi_0(a),\ldots,\xi_{N-1}(a)\right)
=-S(a,x_a).
$$
As $\tilde{\theta}_i(t)=\tilde{u}(t+hi)$, 
$i=0,\ldots,N-1$ and $t\in[a,a+h]$, then 
$$
\tilde{u}(t)
=u^*\left(t,\tilde{x}(t),\tilde{x}(t-r),\partial_2S(t,\tilde{x}(t))\right),
$$
$t\in[a,b]$, is an optimal control that lead us to the minimal cost
$$
C_D[\tilde{u}]=-S(a,x_a).
$$ 
This completes the proof of Theorem~\ref{theo_suf_NLD}.

% -------------------------------------------

\section{Illustrative example}
\label{sec:5}

Let us consider the following problem
studied by G\"{o}llmann et al. in \cite{Gollmann}:
\begin{equation}
\label{example1}
\begin{split}
\min\ \ \ & C[u]=\int_{0}^{3} \left[x^2(t)+u^2(t)\right]dt,\\
\text{s.t.}\ \ & \dot{x}(t)=x(t-1)u(t-2),\\
& x(t)=1,\quad t\in[-1,0],\\
& u(t)=0,\quad t\in[-2,0[,
\end{split}
\end{equation}
which is a particular case of our non-linear delayed optimal control problem $(NLD)$
with $n=m=1$, $a=0$, $b=3$, $r=1$, $s=2$, $g^0(x(3))=0$, $f^0(t,x,y,u,v)=x^2+u^2$
and $f(t,x,y,u,v)=yv$. In \cite{Gollmann}, necessary optimality conditions
were proved and applied to \eqref{example1}. The following candidate 
$(x^*(\cdot),u^*(\cdot))$ was found:
\begin{equation}
\label{eq:ex:cand:x}
\begin{split}
x^*(t)=
\begin{cases}
1,\ & t\in[-1,2],\\[0.8em]
\displaystyle\frac{e^{t-2}+e^{4-t}}{e^2+1},\ & t\in\ [2,3],
\end{cases}
\end{split}
\end{equation}
and
\begin{equation}
\label{eq:ex:cand:u}
\begin{split}
u^*(t)=
\begin{cases}
0,\ & t\in[-2,0[,\\[0.8em]
\displaystyle\frac{e^t-e^{2-t}}{e^2+1},\ & t\in[0,1],\\[0.8em]
0,\ & t\in\ [1,3].
\end{cases}
\end{split}
\end{equation}
It remains missing in \cite{Gollmann}, however, a proof that such candidate 
\eqref{eq:ex:cand:x}--\eqref{eq:ex:cand:u}
is a solution to the problem. It follows from our
sufficient optimality condition that such claim is indeed true. 
In particular, direct computations show that:

\begin{proposition}
\label{eq:prop1}
Function
\begin{equation*}
\begin{split}
S(t,x)=
\begin{cases}
\eta_1(t)x+c_1(t),\ & t\in[0,1],\\[0.8em]
\eta_2(t)x+c_2(t),\ & t\in[1,2],\\[0.8em]
\eta_3(t)x+c_3(t),\ & t\in[2,3],
\end{cases}
\end{split}
\end{equation*}
with
\begin{equation*}
\begin{split}
\eta_1(t) &= -2t+5+\displaystyle\frac{2(e^2-1)}{(e^2+1)^2},\\
\eta_2(t) &= -\left(\displaystyle\frac{4e^2}{(e^2+1)^2}
+2\right)t+\displaystyle\frac{4(e^2-1)}{(e^2+1)^2}
+6+\displaystyle\frac{e^{2t-2}-e^{6-2t}}{(e^2+1)^2},\\
\eta_3(t) &= \displaystyle\frac{2(e^{4-t}-e^{t-2})}{e^2+1},
\end{split}
\end{equation*}
and
\begin{equation*}
\begin{split}
c_1(t)&=\frac{2t(3e^4+4e^2+3)+e^{2t}-e^{4-2t}-15e^4-32e^2-9}{2(e^2+1)^2},\\
c_2(t)&=\frac{2t(3e^4+10e^2+3)+2(e^{6-2t}-e^{2t-2})-17e^4-44e^2-7}{2(e^2+1)^2},\\
c_3(t)&=\frac{4e^2(t-3)+5(e^{2t-4}-e^{8-2t})}{2(e^2+1)^2},
\end{split}
\end{equation*}
is solution of the Hamilton--Jacobi equation \eqref{eq_HJ_main_theo} with $S(3,x^*(3))=0$.
\end{proposition}

% -------------------------------------------
 
\section*{Acknowledgements}

This work is part of first author's Ph.D. project, 
which is carried out at the University of Aveiro
under the support of the Portuguese Foundation 
for Science and Technology (FCT),
fellowship PD/BD/114184/2016.
It was also supported by FCT
within projects UID/MAT/04106/2013 (CIDMA) 
and PTDC/EEI-AUT/2933/2014 (TOCCATA), funded by Project 
3599 -- Promover a Produ\c{c}\~ao Cient\'{\i}fica e Desenvolvimento
Tecnol\'ogico e a Cons\-ti\-tu\-i\-\c{c}\~ao de Redes Tem\'aticas 
and FEDER funds through COMPETE 2020, Programa Operacional
Competitividade e Internacionaliza\c{c}\~ao (POCI).
Silva is also supported by the FCT post-doc 
fellowship SFRH/BPD/72061/2010. 

The authors are grateful to an anonymous referee
for a thorough review, and for suggestions
that significantly contributed to improve the paper.

% -------------------------------------------

% ----------------------------------------------

\end{document}